\documentclass[letterpaper, 10 pt, conference]{ieeeconf}
\IEEEoverridecommandlockouts                              

\usepackage[utf8]{inputenc}
\usepackage{cite}
\usepackage{graphicx} 
\usepackage{amsmath}
\usepackage{siunitx}
\usepackage{multirow}
\usepackage{amssymb}
\usepackage{comment} 
\usepackage{authblk}
\usepackage{algorithmicx}
\usepackage[ruled]{algorithm2e}
\usepackage[table]{xcolor}
\usepackage{lettrine}
\usepackage{url}
 \usepackage{nopageno}

\allowdisplaybreaks[4]
\newcommand\mcedit[1]{\textcolor{black}{#1}}
\newcommand\mdsedit[1]{\textcolor{black}{#1}}

\title{\LARGE \bf
Robust trajectory optimisation for transitions of tiltwing VTOL aircraft}

\author{Martin Doff-Sotta$^{\star}$\qquad Mark Cannon$^{\star}$\qquad  Marko Bacic$^{\star, \dagger}$
\thanks{$^{\star}$All authors are with the Control Group, University of Oxford, Parks Road, Oxford OX1 3PJ, United Kingdom. Corresponding author: Martin Doff-Sotta ({\tt\small martin.doff-sotta@eng.ox.ac.uk})}%
\thanks{$^{\dagger}$On part-time secondment from Rolls-Royce plc.}%
}

\begin{document}

\maketitle
\pagestyle{plain}

\begin{abstract}
We propose a method to generate robust and optimal trajectories for the transition of a tiltwing Vertical Take-Off and Landing (VTOL) aircraft leveraging concepts from convex optimisation, tube-based nonlinear Model Predictive Control (MPC) and Difference of Convex (DC) functions decomposition. 
The approach relies on computing DC decompositions of dynamic models in order to exploit convexity properties and develop a tractable robust optimisation that solves a sequence of convex programs converging to a local optimum of the trajectory generation problem.
%
%
The algorithm developed is applied to an Urban Air Mobility case study. The resulting solutions are robust to approximation errors in dynamic models and provide safe trajectories for aggressive transition manoeuvres at constant altitude.
\\

\end{abstract}

\section{Introduction}
 This paper presents a robust MPC methodology for the trajectory optimisation of VTOL aircraft. Although we consider here the problem of tilt-wing aircraft transition, the method described is equally applicable to tilt-rotors and other forms of VTOL aircraft. 
 
 One of the main challenges associated with VTOL aircraft is stability and control during transition between powered lift and wing-borne flight. This can be problematic as the aircraft experiences large changes in the effective  angle of attack during such manoeuvres.  Achieving successful transitions requires robust flight control laws along feasible trajectories. The computation of the flight transition trajectory is a difficult NonLinear Program (NLP) as it involves nonlinear flight dynamics.



Several attempts were proposed to solve this problem For example, in \cite{UMich19}, the trajectory optimisation for take-off is formulated as a constrained optimisation problem \mcedit{and} solved using NASA's OpenMDAO framework and the SNOPT gradient-based optimiser. The problem of determining minimum energy speed profiles for the forward transition manoeuvre of the Airbus A$^3$ Vahana was addressed in \cite{Iowa19}, considering various phases of flight (cruise, transition, descent). Forward and backward optimal transition manoeuvres at constant altitude are computed in \cite{leo} for a tiltwing aircraft, considering leading-edge fluid injection active flow control and the use of a high-drag device. The main drawback of these approaches is the computational burden associated with solving a NLP, which makes them unsuitable for real-time implementation. 

Another strategy to compute the transition relies on linearisation and convex optimisation resulting in approximate but computationally tractable algorithms. In \cite{ACC}, the transition for a tiltwing VTOL aircraft was computed using convex optimisation by introducing a small angle approximation. This provides a computationally efficient optimisation that could potentially be leveraged online, e.g. for collision avoidance or MPC. The obvious limitation of the approach is the assumption of small angles of attack, which restricts considerably the type of achievable manoeuvres.

While the method of \cite{ACC} introduces a linearisation of the dynamics, there is no consideration for the effect of linearisation error on the dynamics .  In this work, we propose a solution to this problem based on a DC decomposition of the nonlinear dynamics. This allows us to obtain tight bounds on the linearisation error and treat this error as a disturbance in a robust optimisation framework, exploiting an idea from tube MPC \cite{DC-TMPC}. The main idea is to successively linearise the dynamics around predicted trajectories and treat the linearisation error as a bounded disturbance. Due to the DC form of the dynamics, the linearised functions are convex, and so are their linearisation errors. These errors can thus be bounded tightly since they take their maximum at the boundary of the domain, and the trajectories of model states can be bounded by a set of convex inequalities (or tubes \cite{kouvaritakis2016model}). These inequalities form the basis of a computationally-tractable convex tube-based optimisation for the trajectory generation of VTOL aircraft.


The contribution of this research is twofold: i) we solve an open problem in trajectory optimisation  of VTOL aircraft by allowing aggressive transitions at high angle of attack while guaranteeing safety and computational tractability of the scheme; ii) we make a connection between DC decomposition  and robust tube based optimisation and demonstrate the applicability and generalisability of the procedure in \cite{DC-TMPC}. 

This paper is organised as follows. We start by developing a mathematical model of a tiltwing VTOL aircraft in Section \ref{sec:modeling}. In Section \ref{sub:convex_formulation}, we formulate the trajectory optimisation problem and discuss a series of simplifications to obtain a convex program, leveraging ideas from DC decomposition and robust tube MPC. Section \ref{sec:results} discusses simulation results obtained for a case study based on the Airbus A$^3$ Vahana. 
\mcedit{Section~\ref{sec:conclusion} presents conclusions.}

\section{Modeling}
\label{sec:modeling}
\begin{figure}[h]
    \centering
    \includegraphics[width=0.45\textwidth, trim={1.5cm 3.5cm 1.5cm 2cm}, clip]{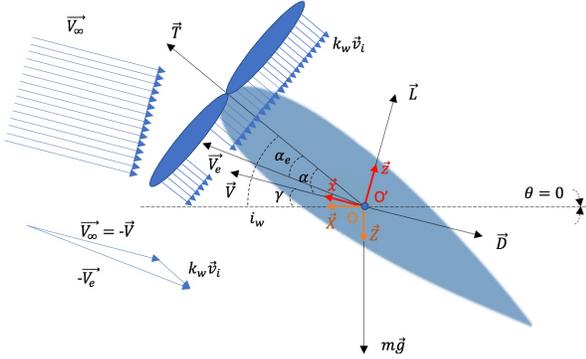} 
    \caption{\mcedit{Force and velocity definitions for a VTOL aircraft}}
    \label{fig:diagram}
\end{figure}

Consider a  longitudinal point-mass model of a tiltwing VTOL aircraft equipped with propellers as shown in Figure~\ref{fig:diagram} and subject to a wind gust disturbance. The Equations Of Motion (EOM) are given in polar form by \cite{ACC}
\begin{alignat}{2}
m \dot{V} &= T \cos{\alpha} - D -mg \sin{\gamma}, 
&\quad 
V(t_0) &= V_{0},
\label{eq:eom1}
\\
m V \dot{\gamma} &= T \sin{\alpha} + L -mg \cos{\gamma}, 
&\quad 
\gamma(t_0) &= \gamma_{0},
\label{eq:eom2}
\end{alignat}
\begin{equation}
J_w \ddot{i}_w = M, \quad i_w(t_0) = i_{0}, \quad \dot{i}_w(t_0) = \Omega_{0},
\label{eq:tiltwing_dyna}
\end{equation}
\begin{equation}
\dot{x} = V \cos{\gamma},
\qquad
\dot{z} = -V \sin{\gamma},
\end{equation}
where the control inputs are the thrust magnitude $T$ and the total torque $M$ delivered by the tilting actuators, and the model states are the aircraft velocity magnitude $V$,  the flight path angle $\gamma$ (defined as the angle of the velocity vector from horizontal), the tiltwing angle $i_w$ and its derivative $\dot{i}_w$, and the position $(x, z)$ with respect to inertial frame $O_{XZ}$. Additional variables are the lift force $L$, drag force $D$ and the angle of attack $\alpha$. All model parameters are defined in Table 1. 

The following input and state constraints apply \cite{ACC}
\begin{gather}
\label{eq:cstr1}
i_w + \theta = \alpha + \gamma, \, \underline{M} \leq M \leq \overline{M}, \\
\label{eq:cstr1bis}
 0 \leq T \leq \overline{T}, \quad 0 \leq V \leq \overline{V},
\\
\label{eq:cstr2}
V(t_0) = V_0
\text{ and }
V(t_f) = V_f,
\\
\label{eq:cstr3}
\underline{a} \leq \dot{V} \leq \overline{a}.
\end{gather}
\mcedit{Here}
$\theta$ is the pitch angle, defined as the angle 
of the fuselage axis from horizontal. For passenger comfort, $\theta$ is regulated via the elevator to track a constant reference $\theta^* = 0$.

In order to account for the effect of the propeller wake on the wing, the flow velocity downstream is augmented by the induced velocity of the propeller. This allows us to define the effective velocity $V_e$ and effective angle of attack $\alpha_e$ seen by the wing as \cite{ACC}
\begin{equation}
 \alpha_e = \arcsin{\left( \frac{V}{V_e} \sin{\alpha}\right)},
\end{equation}
\begin{equation}
V_e = \sqrt{V^2 + \frac{2T}{\rho A n}}. 
\end{equation}

Assuming that the wing is fully immersed in the wake, and that $\alpha_e\ll 1$ to avoid operating the wing in dangerous near-stall regimes\footnote{This will be imposed through a constraint in the optimisation and will be verified \textit{a posteriori} from simulation results.}, the lift and drag are modeled as follows~\cite{ACC}
\begin{equation}
D = \tfrac{1}{2}  \rho S(a_2 \alpha_e^2 + a_1 \alpha_e  + a_0) V_e^2 \approx \tfrac{1}{2}  \rho S(a_1 \alpha_e  + a_0) V_e^2, 
\end{equation}
\begin{equation}
L = \tfrac{1}{2}  \rho S(b_1 \alpha_e  + b_0) V_e^2,
\label{eq:L}
\end{equation}
where $S$ is the wing area, $\rho$ is the air density, $a_0, a_1, a_2$ and $b_0, b_1$ are constant parameters.

\section{Convex optimisation}
\label{sub:convex_formulation}

This paper considers how to robustly generate minimum power trajectories for the transition between powered lift and cruise flight modes, suggesting the following objective function
\begin{equation}
    \label{eq:obj}
    J = \int_{t_0}^{t_f} P/\overline{P}\ \mathrm{d}t,  
\end{equation}
where $P=TV \cos \alpha$ is the drive power and $\overline{P} = \overline{T} \overline{V}$. The optimisation problem consists of minimising \eqref{eq:obj} while satisfying dynamical constraints, input and state constraints \eqref{eq:eom1}-\eqref{eq:L}. As such, this problem is a NLP and we thus consider below how to reformulate the problem as a sequence of convex programs. We introduce 4 key manipulations to do so: i) assuming that a path is known \textit{a priori}, we introduce a change of differential operator to integrate the EOM over space, thus simplifying the structure of the problem; ii) to reduce the couplings between the optimisation variables, we combine both EOM to separate the optimisation of the velocity and torque from the other variables, allowing us to solve 2 smaller optimisation problems sequentially and accelerate computation; iii) we discretise the problem; iv) we approximate the nonlinear dynamics by a difference of convex functions and exploit the fact that convex functions can be bounded tightly by a combination of convex and linear bounds. 

\subsection{Change of differential operator}

Assuming that a path $(x(s), z(s))$ parameterised by the curvilinear abscissa $s$ is known \textit{a priori} (which is usually the case in a UAM context where flight corridors are prescribed) and applying the change of differential operator \cite{bobrow} $\frac{\mathrm{d}}{\mathrm{d}t} = V \frac{\mathrm{d}}{\mathrm{d}s}, \forall V \neq 0$, the dynamics in \eqref{eq:eom1}-\eqref{eq:tiltwing_dyna} can be reformulated as
\small
\begin{align}
\label{eq:eom7}
\frac{1}{2}m E' = T \cos{\alpha} - \frac{1}{2}\rho S\left(a_1 \alpha_e + a_0\right)\left(E + \frac{2T}{\rho A n}\right) -mg \sin{\gamma^*},
\\
\label{eq:eom8}
m E \gamma^{*\prime} = T \sin{\alpha}  + \frac{1}{2}\rho S\left(b_1 \alpha_e + b_0\right)\left(E + \frac{2T}{\rho A n}\right) -mg \cos{\gamma^*}, \\
\label{eq:d2i_w2}
J_w ( \frac{1}{2}E' i_w' +  E i_w'')  = M, \,  i_w(s_0) = i_{0}, \,  i_w'(s_0) \sqrt{E(s_0)}  = \Omega_{0}, 
\end{align}
\normalsize
where $\tfrac{\mathrm{d}\,\cdot}{\mathrm{d}s}=\cdot'$ and $E=V^2$. The flight path angle $\gamma^* = \arctan{(-\mathrm{d}z/\mathrm{d}x)}$ is known \textit{a priori} from the path. 

\subsection{Problem separation}
We next reduce the couplings between the states and inputs in the EOM \eqref{eq:eom7}-\eqref{eq:d2i_w2} by eliminating the angle of attack from the formulation and separating the optimisation into two  subproblems as follows. 
\mcedit{Let} $\lambda = a_1/b_1$, \mcedit{then}  
the combination
\eqref{eq:eom7} $+$ $\lambda$\eqref{eq:eom8} 
yields
\begin{multline}
\!\!\!\!\! \frac{1}{2}m E' +  \underbrace{(\lambda m  \gamma^{*\prime} + \frac{1}{2}\rho S (a_0 -\lambda b_0))}_{c(\gamma^{*\prime})}E  + \underbrace{mg( \sin{\gamma^*} + \lambda \cos{\gamma^*})}_{d(\gamma^*)} 
\\
= \underbrace{T \cos{\alpha} + \lambda T \sin{\alpha} - S^\star (a_0 - \lambda b_0 ) T}_{\tau},
\label{eq:nocvx}
\end{multline}
where 
$S^\star=\frac{S}{An}$ 
\mcedit{and $\tau$  is a virtual input defined by}
\begin{equation}
\tau = T \cos{\alpha} + \lambda T \sin{\alpha} - S^\star (a_0 - \lambda b_0 ) T.
\label{eq:tau}
\end{equation}
The state and input constraints in \eqref{eq:cstr1bis}-\eqref{eq:cstr3} can be rewritten as 
\begin{gather}
\label{eq:cstr_cvx1}
0 \leq E \leq \overline{V}^2, \quad \underline{a} \leq E'/2 \leq \overline{a}, \quad 0\leq \tau \leq \overline{T},
\\
E(s_0) = V_{0}^2
\text{ and }
E(s_f) = V_{f}^2.
\label{eq:cstr_cvx2}
\end{gather}
In the thrust constraint in \eqref{eq:cstr_cvx1}, $\tau$  was chosen as a proxy for $T$ since \mdsedit{$\lambda \ll 1$, and $S^\star (a_0 - \lambda b_0 ) \ll 1$, implying $\tau \approx T \cos \alpha $. This results in the constraint $\tau\leq \overline{T}$ being a relaxed version of the original (we note that the original thrust constraint is inactive in practice -- see Section 5)}. Likewise, the minimum power criterion in \eqref{eq:obj} can be approximated by a convex objective function under these conditions.  By the change of differential operator we obtain \mdsedit{
\begin{equation}
\label{eq:obj_cvx}
    J = \int_{s_0}^{s_f} \tau/\overline{P} \ \mathrm{d}s.
\end{equation}}
Since $\gamma$ and $\gamma'$ are prescribed \mcedit{by} the path, \eqref{eq:nocvx} is a linear equality constraint and the following convex optimisation problem can be constructed to minimise \eqref{eq:obj_cvx} subject to \eqref{eq:nocvx}, \eqref{eq:cstr_cvx1} and \eqref{eq:cstr_cvx2} as follows
\begin{align*}
& \mathcal{P}_1 : \min_{\substack{\tau,\,E,\,a}} & & \int_{s_0}^{s_f} \tau/\overline{P} \, \mathrm{d}s, 
\\
&  \text{ s.t.} & & \frac{1}{2}m E' + c(\gamma^{*\prime}) E + d(\gamma^*) = \tau, 
\\
& & &    0 \leq \tau \leq \overline{T},\   \underline{a} \leq \frac{1}{2} E' \leq \overline{a}, 
\\
& & &   0 \leq E \leq \overline{V}^2,\  E(s_0) = V_{0}^2,\   E(s_f) = V_{f}^2.
\end{align*}
 
Solving $\mathcal{P}_1$ yields the optimal velocity profile
along the path and provides a proxy for the optimal thrust.
However, a
tiltwing angle profile that meets
the dynamical constraints and follows the desired path with $\gamma \approx \gamma^*$ \mcedit{must also be computed}.  
\mcedit{To achieve this}
we use the solution of $\mathcal{P}_1$ to define
a new optimisation problem with variables $\gamma$, $\alpha$, $i_w$,  and $M$ satisfying the constraints \eqref{eq:cstr1}, \eqref{eq:d2i_w2} and, using \eqref{eq:tau} to eliminate the thrust in \eqref{eq:eom8},
%
\begin{align}
m E \gamma' &= \tau \sin{\alpha} -mg \cos{\gamma}
\nonumber
\\
&\quad +\! \frac{1}{2}\rho S \biggl[ b_1 \arcsin \biggl( \frac{\sqrt{E} \sin{\alpha}}{\sqrt{E + \smash{\frac{2\tau}{\rho A n}}\rule{0pt}{8.5pt}}}\biggr) \!+\! b_0\biggr] \Bigl(E \!+\! \frac{2\tau}{\rho A n}\Bigr), 
\nonumber
\\ 
& = f(\alpha, E, \tau) -mg \cos{\gamma}, \label{eq:eom_new}
\end{align}
\mcedit{in which the objective is to minimise the cost function}
\begin{equation}
\label{eq:obj2}
    J_{\gamma} = \int_{s_0}^{s_f} \frac{(\gamma - \gamma^*)^2}{\sqrt{E}} \ \mathrm{d}s. 
\end{equation}
Note that only the two EOM \eqref{eq:eom8} and \eqref{eq:d2i_w2} are needed to construct this new problem since the linear combination \eqref{eq:eom7} + $\lambda$ \eqref{eq:eom8} is enforced with $\tau$ and $E$ prescribed from problem $\mathcal{P}_1$.  We thus state the following optimisation problem
\begin{equation}
\begin{aligned}
& \mathcal{P}_2: \min_{\substack{\alpha,\,\gamma,\,i_w,\, M}} & & \int_{s_0}^{s_f} \frac{(\gamma - \gamma^*)^2}{\sqrt{E}} \ \mathrm{d}s
\nonumber\\
&  \text{ s.t.} & & m E \gamma' =f(\alpha, E, \tau) -mg \cos{\gamma},  
\nonumber\\
& & &  J_w (\tfrac{1}{2} E' i_w' + i_w'' E)  = M,\ i_w(s_0) = i_{0},
\nonumber
\\
& & &   i_w'(s_0) \sqrt{E(s_0)}  = \Omega_{0},
\nonumber\\
& & &    i_w  = \alpha + \gamma, 
\nonumber\\
& & &    \underline{M} \leq M \leq \overline{M},\   \underline{\alpha} \leq \alpha \leq \overline{\alpha}
\nonumber\\
& & &   \underline{\gamma} \leq \gamma \leq \overline{\gamma},\  \underline{i_w}  \leq i_w \leq \overline{i_w},
\nonumber\\
\end{aligned}
\end{equation}
and reconstruct the input $T$ and state $V$ \textit{a posteriori} using \eqref{eq:tau} and $V = \sqrt{E}$.  Given the solution of both problems as functions of the independent variable $s$, the \mcedit{final} step is to map the solution to time domain by reversing the change of differential operator and integrating
\[
t(\xi) = \int_{s_0}^\xi \frac{\mathrm{d}s }{V(s)}.
\]


We have now achieved the separation into two subproblems $\mathcal{P}_1$ and $\mathcal{P}_2$, as described in \cite{ACC}. 

\subsection{Discretisation}
\label{sec:discrete}
The decision variables in $\mathcal{P}_1$ and $\mathcal{P}_2$ are functions defined on the interval $[s_0,s_f]$. To obtain computationally tractable problems, we consider $N+1$ discretisation points $\{s_0, s_1,  \ldots, s_{N}\}$ of the path, with spacing $\delta_k = s_{k+1} - s_{k}$, $k = 0, \ldots, N-1$ ($N$ steps). The notation $\{u_0, \ldots, u_{N}\}$ is used for the sequence of the discrete values of a continuous variable $u$ evaluated at the discretisation points of the mesh, where $u_k = u(s_{k})$, $\forall k \in \{0, \ldots, N\}$.

Assuming a path $s_k \rightarrow (x_k, z_k)$, the prescribed flight path angle and rate are discretised as follows
\begin{align}
\gamma_k^* &= \arctan \Bigl(-\frac{z_{k+1}- z_{k}}{x_{k+1}- x_{k}}\Bigr), \quad k \in \{0, \ldots, N-1\},
\label{eq:gamma}
\\
\gamma_k^{*\prime} &= \begin{cases} (\gamma_{k+1}^*- \gamma_{k}^*)/{\delta_k}, &  k \in \{0, \ldots, N-2\}, \\
\gamma_{N-2}^{*\prime}, & k = N-1.
\end{cases}
\label{eq:dgamma}
\end{align}
The resulting discretised versions of $\mathcal{P}_1$ and $\mathcal{P}_2$ are
\begin{equation}
\begin{aligned}
& \mathcal{P}_1^\dagger : \min_{\substack{\tau,\,E,\,a}} & & \sum_{k=0}^{N-1} \tau_k/\overline{P} \, \delta_k,
\nonumber\\
&  \text{ s.t.} & & E_{k+1}  = E_{k} + \frac{2 \delta_{k}}{m} (\tau_k - c(\gamma_k^{*\prime}) E_k - d(\gamma_k^*)), 
\nonumber\\
& & &    0 \leq \tau_k \leq \overline{T}, \quad  \underline{a} \leq \frac{E_{k+1} - E_{k}}{2 \delta_{k}} \leq \overline{a}, 
\nonumber\\
& & &   0 \leq E_k \leq \overline{V}^2, \quad  E_0 = V_{0}^2, \quad   E_{N} = V_{f}^2,
\nonumber\\
\end{aligned}
\end{equation}
\begin{equation}
\begin{aligned}
& \mathcal{P}_2^\dagger : \min_{\substack{\alpha,\,\gamma,\\i_w,\, \zeta,\, M}} & & \sum_{k=0}^{N-1} \frac{(\gamma_k - \gamma_k^*)^2}{\sqrt{E_k}}\delta_k 
\nonumber\\
& \text{s.t.} & &   \gamma_{k+1} = \gamma_{k} + \frac{\delta_k}{m E} (f_k(\alpha_k, E_k, \tau_k) -f_{\gamma_k}(\gamma_k)), 
\nonumber\\
& & &    i_{w,k}  = \alpha_k + \gamma_k, 
\nonumber\\
& & & i_{w,k+1}   = i_{w,k} + \zeta_k \delta_k, \quad i_{w,0} = i_{0}, 
\nonumber\\
& & & \zeta_{k+1}  = \zeta_{k}\Bigl(1-\frac{E_{k+1} - E_k}{2E_k}\Bigr) + \frac{M_k \delta_k}{J_w {E_k}}, 
\nonumber \\
& & & \zeta_0 \sqrt{E}_0  = \Omega_{0},
\nonumber \\
& & &    \underline{M} \leq M_k \leq \overline{M}, \quad \underline{\alpha} \leq \alpha_k \leq \overline{\alpha}, 
\nonumber\\
& & &   \underline{\gamma} \leq \gamma_k \leq \overline{\gamma}, \quad   \underline{i_w}  \leq i_{w,k} \leq \overline{i_w}.
\nonumber\\
\end{aligned}
\end{equation}
%
where $f_{\gamma_k} = -mg \cos{\gamma_k}$. The input and state \mcedit{variables} are reconstructed using
\begin{equation}
  \label{eq:reconstruct1}
    T_k = \frac{\tau_k}{ \cos{\alpha_k} + \lambda  \sin{\alpha_k} - S^\star (a_0 - \lambda b_0 )},  \quad V_k = \sqrt{E_k} ,  
\end{equation}
\mcedit{and} 
the time $t_k$ associated with each discretisation point 
\mcedit{is computed, allowing solutions to be expressed}
as time series
\begin{equation}
\label{eq:back2time}
t_k = \sum_{j=0}^{k-1} \frac{\delta_j}{V_j}.
\end{equation}

We now have a pair of finite dimensional problems $\mathcal{P}_1^\dagger$ and $\mathcal{P}_2^\dagger$, but the latter is still nonconvex due to the nonlinear functions $f_k$ and $f_{\gamma_k}$ in the dynamics. On a restricted domain $\gamma_k \in [-\pi/2; \pi/2]$, $f_{\gamma_k}$ is a convex function of $\gamma_k$, making it possible to derive tight convex bounds on $f_{\gamma_k}$ (as discussed in Section \ref{sec:cvx_relax}). However, this is not the case with $f_k$ and we introduce a method to alleviate this limitation in what follows.  

\subsection{DC decomposition}

Motivated by the fact that convex functions can be bounded tightly by convex and linear inequalities (as in \cite{DC-TMPC}), we seek a decomposition of $f$ as a Difference of Convex (DC) functions: $f_k = g_k -h_k$, where $g_k, h_k$ are convex. A DC decomposition always exists if $f_k \in \mathcal{C}^2$ \cite{hartman}. 

Note that since $(E_k, \tau_k)$ are obtained from problem $\mathcal{P}_1^\dagger$, the function $f_k(\alpha_k, E_k, \tau_k)$ is single-valued (in $\alpha_k$) which considerably simplifies the task of finding a DC decomposition, and motivates the above approach of separating the initial problem in two subproblems with fewer couplings between the variables. However, $f_k$ is also time varying through its dependence on parameters $(E_k, \tau_k)$ generated online. This requires us to find a DC split for every instance of $(E_k, \tau_k), \, \forall k \in [0, 1, ..., N]$ which can be intractable if the horizon is large or the sampling interval is small.  

Instead, we adopt the more pragmatic approach of i) precomputing offline the DC decompositions on a downsampled grid of values $(E_i, \tau_j), \, \forall (i, j) \in [0, 1, ..., N_s] \times [0, 1, ..., M_s]$ where $N_s, M_s \ll N+1$ and ii) interpolating the obtained decompositions online using a lookup table. 

\subsubsection{Precomputation of the DC decomposition} Inspired by \cite{ahmadi}, we develop a computationally tractable method for the DC decomposition of a function $f_k(\alpha)$ based on an approximation\footnote{Note that any continuous function can be approximated arbitrarily closely by a polynomial.} of the function by a polynomial of degree $2n$: 
\begin{equation}
\label{eq:pol_approx}
f_k (\alpha) \approx p_{k, 2n} \alpha^{2n} + ... + p_{k, 1} \alpha + p_{k, 0} = y^\top P_k y, 
\end{equation}
where $y = [1 \quad \alpha \quad ... \quad \alpha^n]^\top \in \mathbb{R}^{n+1}$ is a vector of monomials of increasing order and $P_k = P_k^\top \in R^{n+1 \times n+1}$ is the Gram matrix of the polynomial defined by $\{P_k\}_{ij} = p_{k, i+j +1}/\lambda(i, j), \, \forall i,j \in [0, 1, ..., n]$ where $\lambda(i, j) = i + j + 1$ if $i+j \leq n $ and $\lambda(i, j) = 2n + 1 - (i + j)$. Given $N_s$ samples $F_{k, s} \, \forall s \in [1, ..., N_s]$ of the function $f_k$,  the polynomial approximation can be obtained by solving a least square problem to find the coefficients that best fit the samples. 

We now seek the symmetric matrices $Q_k$, $R_k$ such that 
\[
y^\top P_k y = y^\top Q_k y - y^\top R_k y,
\]
where $g_k \approx y^\top Q_k y$ and $h_k \approx y^\top R_k y$ are convex polynomials in $\alpha$. Such conditions can be satisfied if the Hessians $d^2g_k/d\alpha^2 = y^\top H_{g_k} y$ and $d^2h_k/d\alpha^2 = y^\top H_{h_k} y$ are Positive Semi-Definite (PSD), i.e. if the following Linear Matrix Inequalities (LMI) hold
\[
H_{g_k} \equiv {D^\top}^2 Q_k + Q_k D^2 + 2 D^\top Q_k D \succeq 0, 
\]
\[
 H_{h_k} \equiv {D^\top}^2 R_k + R_k D^2 + 2 D^\top R_k D \succeq 0, 
\]
where $D$ is a matrix of coefficients such that $dy/d\alpha = D y $. Finding the DC decomposition thus reduces to solving the following Semi Definite Program (SDP)
\begin{equation}
\begin{aligned}
 & \mathcal{SDP} : \min_{\substack{H_{g_k}}} & & \mathrm{tr} \, H_{g_k}
\nonumber\\
&  \text{ s.t.} & & H_{g_k} \succeq 0 ,  
\nonumber\\
& & &  H_{g_k} - ({D^\top}^2 P_k + P_k D^2 + 2D^\top P_k D) \succeq 0 ,
\nonumber
\end{aligned}
\end{equation}
 and computing $H_{h_k} = H_{g_k} - {D^\top}^2 P_k + P_k D^2 + D^\top P_k D$, followed by the double integration $d^2g_k/d\alpha^2 = y^\top H_{g_k} y $ and $d^2h_k/d\alpha^2 = y^\top H_{h_k} y$ to recover $g_k$ and $h_k$. This operation is repeated at each point $(E_i, \tau_j)$ of the grid to assemble a look-up table of polynomial coefficients.  Note that the objective was chosen so as to regularise the solutions for $g_k$, $h_k$ by minimising a proxy for their average curvature, in order to minimise linearisation errors later on. 

\begin{figure}
         \centering
         \includegraphics[width=0.4\textwidth, trim={0cm 5cm 0cm 5cm}]{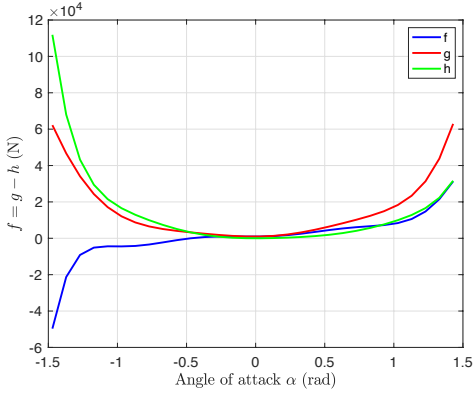}
         \caption{Example of a DC decomposition for a given $k$.}
         \label{fig:split}
\end{figure}

In Figure \ref{fig:split}, we illustrate a typical DC decomposition of the nonlinear dynamics for a given $(E_i, \tau_j)$.

 \subsubsection{Coefficient interpolation} A bilinear interpolation of the coefficients is performed online to obtain the DC decomposition for each $(E_k, \tau_k), \, \forall k \in [0, ..., N]$.  
 This operation preserves convexity since the interpolated polynomial coefficients are a weighted sum of the coefficients in the lookup table. 
\subsection{Convex relaxation}
\label{sec:cvx_relax}
Consider again the nonlinear dynamics in problem $\mathcal{P}_2^\dagger$, using the DC decomposition of $f_k$ computed in the previous section and eliminating the angle of attack via $\alpha_k = i_{w, k} - \gamma_k$ to reduce the number of states, we obtain 
\begin{equation}
\label{eq:dyna_eq_cvx}
\begin{split}
    \gamma_{k+1} = \gamma_{k} + \frac{\delta_k}{m E} (g_k(i_{w, k} - \gamma_k, E_k, \tau_k) \\ - h_k(i_{w, k} - \gamma_k, E_k, \tau_k) -  mg \cos \gamma_k). 
\end{split}
\end{equation}

All nonlinearities in equation \eqref{eq:dyna_eq_cvx} above involve convex and concave functions of the states $i_{w, k}$ and $\gamma_k$ whose dynamics are given by
\begin{gather}
\label{eq:lqr1}
i_{w,k+1}   =  i_{w,k} + \zeta_k \delta_k, \\
\label{eq:lqr2}
\zeta_{k+1}  = \zeta_{k}\Bigl(1-\frac{E_{k+1} - E_k}{2E_k}\Bigr) + \frac{M_k \delta_k}{J_w {E_k}}.
\end{gather}

In what follows we will exploit the convexity properties of the functions $g_k, h_k, f_{\gamma_k} = -mg \cos \gamma_k$ in \eqref{eq:dyna_eq_cvx} to approximate the dynamics by a set of convex inequalities with tight bounds on the state trajectories. To do so, we linearise the dynamics successively around feasible guessed trajectories and treat the linearisation error as a bounded disturbance  \cite{DC-TMPC}. We use the fact that the linearisation error of a convex (resp. concave) function is also convex (resp. concave) and can thus be bounded tightly since its maximum (resp. minimum) occurs at the boundary of the set on which the function is constrained. This allows us to construct a robust optimisation using the tube-based MPC framework \cite{kouvaritakis2016model}, and to obtain solutions that are robust to the model error introduced by the linearisation. 

We start by assuming the existence of a set of feasible trajectories $i_{w, k}$ and $\gamma_k$ for \eqref{eq:dyna_eq_cvx}-\eqref{eq:lqr2} and consider the perturbed dynamics
\begin{equation}
\label{eq:dyna_eq_cvx_perturbed}
\begin{split}
    \gamma_{k+1} = \gamma_{k} + \frac{\delta_k}{m E} (g_k^\circ + \nabla g_k^\circ (i_{w, k} - \gamma_k -(i_{w, k}^\circ - \gamma_k^\circ)) \\ + w_1 - h_k^\circ - \nabla h_k^\circ (i_{w, k} - \gamma_k -(i_{w, k}^\circ - \gamma_k^\circ)) - w_2 \\ -  mg \cos \gamma_k^\circ + mg \sin \gamma_k^\circ (\gamma_k - \gamma_k^\circ) + w_3) . 
\end{split}
\end{equation}
where $g_{k}^{\circ} = g_k(i_{w, k}^\circ - \gamma_k^\circ)$, $h_k^\circ = h_k(i_{w, k}^\circ - \gamma_k^\circ)$ are the functions $g_k, h_k$ evaluated along the guessed trajectory, $\nabla g_k^\circ = dg_k/d\alpha_k (i_{w, k}^\circ - \gamma_k^\circ)$, $\nabla h_k^\circ = dh_k/d\alpha_k (i_{w, k}^\circ - \gamma_k^\circ)$ are the first order derivatives of $g_k, h_k$ evaluated along the guessed trajectory, and $w_1(i_{w} - \gamma, i_{w, k}^\circ - \gamma_k^\circ)$, $w_2(i_{w} - \gamma, i_{w, k}^\circ - \gamma_k^\circ)$, $w_3(\gamma, \gamma_k^\circ)$ are the convex linearisation errors of $g_k, h_k, f_{\gamma_k} = -mg \cos \gamma_k$ respectively. Since these linearisation errors are convex, they take their maximum on the boundary of the set over which the functions are constrained. Moreover, by definition, their minimum on this set is zero (Jacobian linearisation). We thus infer the following relationships $\forall i = \{1, 2\}$ and noting $f_1 \equiv g, f_2 \equiv h$
\begin{equation}
\label{eq:wi_min}
\min_{\substack{\gamma \in [\underline{\gamma}_k, \overline{\gamma}_k]\\ i_w \in [\underline{i}_{w, k}, \overline{i}_{w, k}]}} w_i(i_w - \gamma, i_{w}^\circ - \gamma^\circ) = 0, 
\end{equation}
\begin{equation}
\label{eq:wi_max}
\begin{split}
\max_{\substack{\gamma \in [\underline{\gamma}_k, \overline{\gamma}_k]\\ i_w \in [\underline{i}_{w, k}, \overline{i}_{w, k}]}} w_i(i_w - \gamma, i_{w}^\circ - \gamma^\circ) = \qquad \qquad \qquad \qquad \qquad\\ \max \{f_{i, k} - f_{i, k}^\circ - \nabla f_{i, k}^\circ (\overline{i}_{w, k} - \underline{\gamma}_k - (i_{w, k}^\circ - \gamma_k^\circ)) ;\\ f_{i, k} - f_{i, k}^\circ - \nabla f_{i, k}^\circ (\underline{i}_{w, k} - \overline{\gamma}_k - (i_{w, k}^\circ - \gamma_k^\circ))\},
\end{split}
\end{equation}
\begin{equation}
\label{eq:w3_max}
\begin{split}
\min_{\substack{\gamma \in [\underline{\gamma}_k, \overline{\gamma}_k]}} w_3(\gamma, \gamma^\circ) = 0 \quad \text{and} \quad \max_{\substack{\gamma \in [\underline{\gamma}_k, \overline{\gamma}_k]}} w_3(\gamma, \gamma^\circ) =\\
\max \{-mg \cos \underline{\gamma}_k + mg \cos \gamma_k^\circ - mg \sin \gamma_k^\circ (\underline{\gamma}_k - \gamma_k^\circ); \\
-mg \cos \overline{\gamma}_k + mg \cos \gamma_k^\circ - mg \sin \gamma_k^\circ (\overline{\gamma}_k - \gamma_k^\circ) \}
\end{split}
\end{equation}
where we assumed that the state trajectories $\gamma_k$ and $i_{w, k}$ lie within "tubes" whose cross-sections are parameterised by means of elementwise bounds $\gamma_{k} \in [\underline{\gamma}_{k}, \overline{\gamma}_{k}]$ and $i_{w, k} \in [\underline{i}_{w, k}, \overline{i}_{w, k}] ,\ \forall k$, which are considered to be optimisation variables. Given these bounds on the states at a given time instant and by virtue of equations \eqref{eq:wi_min}-\eqref{eq:w3_max}, the bounds on the states at the next time instant satisfy the following convex inequalities 
\begin{equation}\label{eq:dyna_ineq5}
\begin{split}
    \overline{\gamma}_{k+1} \geq \max_{\substack{\gamma \in \{\underline{\gamma}_k; \overline{\gamma}_k\}\\ i_w \in \{\underline{i}_{w, k}; \overline{i}_{w, k}\}}} \Big\{ \gamma + \frac{\delta_k}{m E} (g_k(i_{w} - \gamma) - h_k(i_{w}^\circ - \gamma^\circ) \\ - \nabla h_k^\circ (i_{w} - \gamma -(i_{w, k}^\circ - \gamma_k^\circ))  -  mg \cos \gamma) \Big\}, 
\end{split}
\end{equation}
\begin{equation}
\label{eq:dyna_ineq6}
\begin{split}
    \underline{\gamma}_{k+1} \leq \min_{\substack{\gamma \in \{\underline{\gamma}_k, \overline{\gamma}_k\}\\ i_w \in \{\underline{i}_{w, k}, \overline{i}_{w, k}\}}} \Big\{ \gamma + \frac{\delta_k}{m E} (g_k(i_{w}^\circ - \gamma^\circ) \\+\nabla g_k^\circ (i_{w} - \gamma- (i_{w, k}^\circ - \gamma_k^\circ))  - h_k(i_{w} - \gamma) \\-  mg \cos \gamma^\circ + mg \sin \gamma^\circ (i_{w} - \gamma)  ) \Big\}, 
\end{split}
\end{equation}
\begin{equation}
\label{eq:iw}
    \overline{i}_{w,k+1}   \geq  \overline{i}_{w,k} + \zeta_k \delta_k, \quad \underline{i}_{w,k+1}   \leq  \underline{i}_{w,k} + \zeta_k \delta_k.
\end{equation}

These conditions involve only minimisations of linear functions and maximisations of convex functions. Note that the functions to optimise no longer need to be evaluated on continuous intervals but at their boundaries $\{\underline{\gamma}_k, \overline{\gamma}_k\}$ and $\{\underline{i}_{w, k}, \overline{i}_{w, k}\}$ which implies that each maximisation and minimisation above reduces to $2^2 = 4$ convex inequalities. Moreover, this number can be reduced to avoid the curse of dimensionality since the coefficients of the linear functions appearing in each maximisation and minimisation are known. Finally, the computational burden can be further reduced  by introducing a low order approximation of the polynomials in \eqref{eq:dyna_ineq5}-\eqref{eq:iw}. This was obtained by computing, before including the constraints in the optimisation, a series of quadratic polynomials to each $g_k$, $h_k$, $\forall k$ that are a best fit around $i_{w, k}^\circ - \gamma_k^\circ$. 

The tube defined by inequalities \eqref{eq:dyna_ineq5}-\eqref{eq:iw} can be used to replace $\mathcal{P}_2^\dagger$ by a sequence of convex programs. Given the solution of $\mathcal{P}_1^\dagger$ and given a set of feasible (suboptimal) trajectories $i_{w, k}^\circ$, $\gamma_k^\circ$ satisfying \eqref{eq:dyna_eq_cvx}-\eqref{eq:lqr2}, the following convex problem is solved sequentially

\begin{equation}
\begin{aligned}
& \mathcal{P}_2^\ddagger : \min_{\substack{\overline{\gamma},\,\underline{\gamma},\, \overline{i}_w,\\ \underline{i}_w,\, \zeta,\, M, \theta}} & & \sum_{k=0}^{N-1} \frac{\theta_k^2}{\sqrt{E_k}}\delta_k 
\nonumber\\
& \text{s.t.} & &   \theta_k \geq |\overline{\gamma}_k - \gamma_k^*|, \quad \theta_k \geq |\underline{\gamma}_k - \gamma_k^*|, 
\nonumber\\
& & & \overline{\gamma}_{k+1} \geq \max_{\substack{\gamma \in \{\underline{\gamma}_k; \overline{\gamma}_k\}\\ i_w \in \{\underline{i}_{w, k}; \overline{i}_{w, k}\}}} \Big\{ \gamma + \frac{\delta_k}{m E} (g_k(i_{w} - \gamma) \nonumber \\
& & & - h_k(i_{w}^\circ - \gamma^\circ)  - \nabla h_k^\circ (i_{w} - \gamma -(i_{w, k}^\circ - \gamma_k^\circ))  \nonumber\\
& & & -  mg \cos \gamma) \Big\},
\nonumber\\
& & &  \underline{\gamma}_{k+1} \leq \min_{\substack{\gamma \in \{\underline{\gamma}_k, \overline{\gamma}_k\}\\ i_w \in \{\underline{i}_{w, k}, \overline{i}_{w, k}\}}} \Big\{ \gamma + \frac{\delta_k}{m E} (g_k(i_{w}^\circ - \gamma^\circ) \nonumber \\
& & & +\nabla g_k^\circ (i_{w} - \gamma -(i_{w, k}^\circ - \gamma_k^\circ))  - h_k(i_{w} - \gamma) 
\nonumber\\  
& & & -  mg \cos \gamma^\circ + mg \sin \gamma^\circ (i_{w} - \gamma)  ) \Big\},
\nonumber\\
& & & \overline{i}_{w,k+1}   \geq  \overline{i}_{w,k} + \zeta_k \delta_k, \quad \overline{i}_{w,0} = i_{0}, 
\nonumber\\
& & & \underline{i}_{w,k+1}   \leq  \underline{i}_{w,k} + \zeta_k \delta_k, \quad \underline{i}_{w,0} = i_{0}, 
\nonumber\\
& & & \zeta_{k+1}  = \zeta_{k}\Bigl(1-\frac{E_{k+1} - E_k}{2E_k}\Bigr) + \frac{M_k \delta_k}{J_w {E_k}}, 
\nonumber \\
& & & \zeta_0 \sqrt{E}_0  = \Omega_{0},
\nonumber \\
& & &    \underline{M} \leq M_k \leq \overline{M}, \quad  \underline{i_w}  \leq \underline{i}_{w,k}, \quad   \overline{i}_{w,k} \leq \overline{i_w}, 
\nonumber\\
& & &   \underline{\alpha} \leq \underline{i}_{w, k} - \overline{\gamma}_k , \quad  \overline{i}_{w, k} - \underline{\gamma}_k \leq \overline{\alpha}, 
\nonumber\\
& & &   \underline{\gamma} \leq  \underline{\gamma}_k , \quad  \overline{\gamma}_k \leq \overline{\gamma}, \quad \overline{\gamma}_{0} = \gamma_{0}, \quad \underline{\gamma}_{0} = \gamma_{0}. 
\nonumber\\
\end{aligned}
\end{equation}

\mdsedit{After each iteration of this problem, the guessed trajectories are updated by passing $M_k$ through the dynamics \eqref{eq:lqr1}-\eqref{eq:lqr2}, and updating
\begin{equation}
\label{eq:update1}
E_{k+1}  = E_k +  \frac{2 \delta_k}{m} (T_k \cos{\alpha_k} - D_k -mg \sin{\gamma_k}), 
\end{equation} 
\begin{equation}
\label{eq:update2}
\gamma_{k+1} = \gamma_{k} + \frac{\delta_k}{m E} (f_k(\alpha_k, E_k, \tau_k) -f_{\gamma_k}(\gamma_k)), 
\end{equation} 
where $D_k = \tfrac{1}{2}  \rho S(a_2 \alpha_e^2 + a_1 \alpha_e  + a_0) V_e^2$ and $T_k$ is obtained using equation \eqref{eq:reconstruct1}}. The process is repeated until $|\overline{\gamma}_k - \underline{\gamma}_k|$ and $|\overline{i}_{w, k} - \underline{i}_{w, k}|$ have converged. Once $\mathcal{P}_1^\dagger$ and $\mathcal{P}_2^\ddagger$ have been solved, we check whether $|\gamma_k^*-\gamma_k| \leq \epsilon$ $\forall k \in \{0, \ldots, N\}$, where $\epsilon$ is a specified tolerance. If this condition is not met ($\mathcal{P}_2^\ddagger$ may admit solutions that allow $\gamma_k$ to differ from the assumed flight path angle $\gamma^*_k$), the problem is reinitialized with the updated flight path angle and rate $\gamma^*_k \gets \gamma_k$, $\gamma^{*\prime}_k \gets \gamma_k'$ and \mcedit{$\mathcal{P}_1^\dagger$ and $\mathcal{P}_2^\ddagger$ are solved again}.
When \mcedit{the solution tolerance is met}
(or 
\mcedit{the} maximum number of iterations \mcedit{is exceeded}) the problem is considered solved and the 
input and state \mcedit{variables} are reconstructed 
using the equations in \eqref{eq:reconstruct1} and the time $t_k$ associated with each discretisation point 
\mcedit{is computed with \eqref{eq:back2time}, allowing solutions to be expressed}
as time series. 
The procedure is summarised in Algorithm 1. 

Remarks on $\mathcal{P}_2^\ddagger$: i) the angle of attack has been eliminated from the formulation; ii) the slack variable $\theta_k$ was introduced to enforce the objective  $\sum_k \max_{\gamma_k \in \{\underline{\gamma}_k;\overline{\gamma}_k\}} (\gamma_k - \gamma^*)^2/\sqrt{E_k}$; iii) to ensure convexity, it is important that $\overline{\gamma} \leq \pi/2$ and $\underline{\gamma} \geq -\pi/2$; 
iv) to improve numerical stability, $M_k$ can be replaced by $M_k + K_{i, k}(i_{w, k} - i_{w, k}^\circ) + K_{\zeta, k} (\zeta_k - \zeta_{k}^\circ)$ with $i_{w, k} \in \{\underline{i}_{w, k}; \overline{i}_{w, k} \}$ where $K_{i, k}$ and $K_{\zeta, k}$ are gains obtained, e.g. by solving a LQR problem for the time varying linear system in equations \eqref{eq:lqr1} and \eqref{eq:lqr2}; v) order reduction was performed on polynomials $g_k$, $h_k$, i.e. quadratic polynomials were fitted to $g_k$, $h_k$ around $i_{w, k}^\circ - \gamma_k^\circ$ for all $k$ by solving a least squares problem before running the optimisation. 



\begin{algorithm}
\caption{Convex trajectory optimisation}
\label{algo}
    Compute the DC decomposition of $f_k$ on grid of points $(E_i, \tau_j)$: approximate $f_k$ by a polynomial \eqref{eq:pol_approx} and solve $\mathcal{SDP}$ at each point to obtain  $g_k$, $h_k$ and build a look-up table of polynomial coefficients. \\
    \mcedit{Compute $\gamma^*$, $\gamma^{*\prime}$ using \eqref{eq:gamma}, \eqref{eq:dgamma}}
    \mcedit{and initialise: $\gamma\gets\gamma^*$, $\gamma'\gets\gamma^{*\prime}$}, $\gamma^* \gets \infty$, $j \gets 0$\\
    \While{$\displaystyle\max_{k \in \{0, \ldots, N\}}|\gamma_k^*-\gamma_k| > \epsilon$ $\&$ $j < $ MaxIters}{
    $\gamma^{*} \gets \gamma$, $\gamma^{*\prime} \gets \gamma'$\\
    Solve problem $\mathcal{P}_1^\dagger$.\\
    Compute the gains $K_{\zeta, k}$ and $K_{i, k}$. \\
    Compute feasible trajectories $i_{w, k}^\circ$ and $\gamma_{k}^\circ$.\\
    $i \gets 0$\\
    \While{$\displaystyle\max_{k \in \{0, \ldots, N\}}|\overline{\gamma}_k -\underline{\gamma}_k| > \epsilon$ $\&$ $\displaystyle\max_{k \in \{0, \ldots, N\}} |\overline{i}_{w, k} -\underline{i}_{w, k}| > \epsilon$ $\&$ $i < $ MaxIters}{

    \For{$k\gets0$ \KwTo $N$}{
    Interpolate the coefficients of the DC polynomials $g_k, h_k$ from the look-up table. \\
    Fit a quadratic model to the interpolated $g_k, h_k$ that is a best fit at $\i_{w, k}^\circ - \gamma_k^\circ$.\\
    }
    Solve problem $\mathcal{P}_2^\ddagger$. \\
    Update guess trajectories with \eqref{eq:lqr1} - \eqref{eq:lqr2} and  \eqref{eq:update1} - \eqref{eq:update2}. \\
    $i \gets  i +1$\\ 
    }
    $j \gets j+1$
    }
    \For{$k\gets0$ \KwTo $N$}{
    $T_k \gets \frac{\tau_k}{ \cos{\alpha_k} + \lambda  \sin{\alpha_k} - S^\star (a_0 - \lambda b_0 )}$ \\
    $V_k \gets \sqrt{E_k}$\\
    $t_k \gets \sum_{j=0}^{k-1} \frac{\delta_j}{V_j}$
    }
\end{algorithm}

\section{Results}
\label{sec:results}
We consider a case study \mcedit{based on}
the Airbus A$^3$ Vahana. The \mcedit{aircraft}
parameters are reported in Table \ref{tab:param}.
We run Algorithm 1 \mcedit{using}
the convex programming software package CVX \cite{cvx} with \mcedit{the solver} Mosek \cite{mosek} to compute the optimal trajectory for 2 different transition manoeuvres, with boundary conditions given in Table \ref{tab:BC}. For the sake of simplicity, and unless otherwise stated, we limit the number of iterations of problem $\mathcal{P}_1^\dagger$ to 1 and of $\mathcal{P}_2^\ddagger$ to 3. The average computation time per iteration of  $\mathcal{P}_2^\ddagger$ was 7.3s.

\begin{figure}
         \centering
         \includegraphics[width=0.5\textwidth, trim={0cm 5cm 0cm 5cm}, clip]{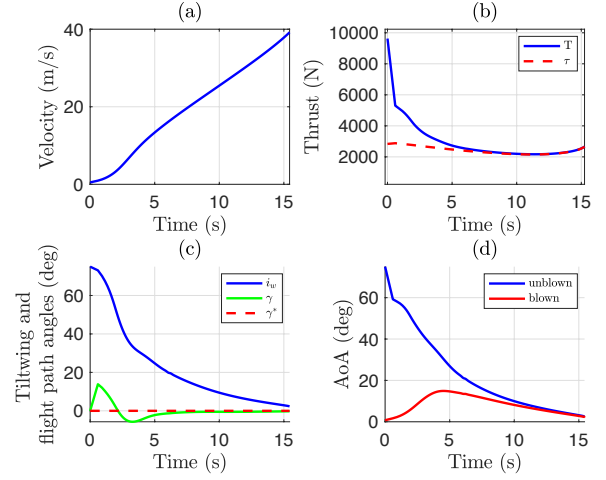}
         \caption{Forward transition (scenario 1). The arrows in the last subplot represent the thrust vector along the trajectory.}
         \label{fig:traj1}
\end{figure}

\begin{figure}
         \centering
         \includegraphics[width=0.4\textwidth, trim={0cm 5cm 0cm 5cm}]{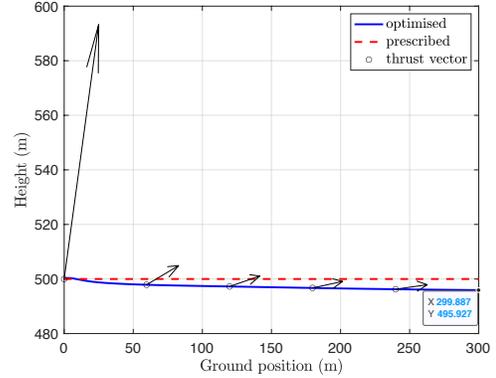}
         \caption{Altitude variation and thrust vector field during the forward transition (scenario 1).}
         \label{fig:zoom}
\end{figure}

The first scenario is a (near) constant altitude forward transition. This manoeuvre is abrupt and requires a zero flight path angle throughout as illustrated in Figure \ref{fig:traj1}. As the aircraft transitions from powered lift to cruise, the velocity magnitude increases (a) and the thrust decreases (b), illustrating the change in lift generation from propellers to wing. The tiltwing angle drops quickly at the beginning (c), resulting in an increase in the angle of attack (d). The slight discrepancy in the flight path angle curves in (c) illustrates that problem $\mathcal{P}_2^\ddagger$ needs not necessarily generate a flight path angle profile corresponding to the exact desired path if the latter is not feasible. Note from graph (d) that the effective angle of attack stays within reasonable bounds, indicating that the wing is not stalled. By contrast to the solution presented in \cite{ACC}, the angle of attack is not constrained to small values and we can thus achieve a more aggressive transition at an almost constant altitude, with a maximum altitude drop of about 4\,m, see Figure \ref{fig:zoom}.    

\begin{figure}
         \centering
         \includegraphics[width=0.4\textwidth, trim={0cm 5cm 0cm 5cm}]{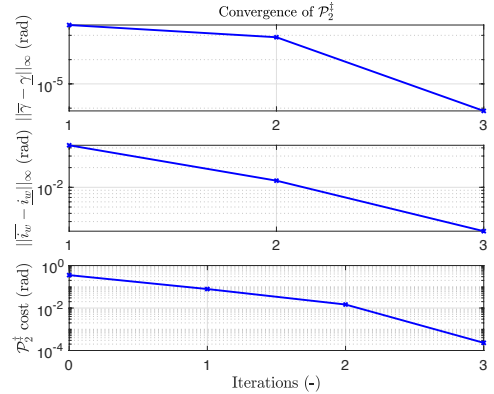}
         \caption{Convergence of tube bounds and objectivefor $\mathcal{P}_2^\ddagger$ (scenario 1).}
         \label{fig:obj}
\end{figure}

Convergence of problem $\mathcal{P}_2^\ddagger$ after 3 iterations is shown in Figure \ref{fig:obj}. The tube bounds and objective converge quickly toward infinitesimal values after just a few iterations. After that, no more progress can be achieved. 

\begin{figure}
         \centering
         \includegraphics[width=0.5\textwidth, trim={0cm 5cm 0cm 5cm}, clip]{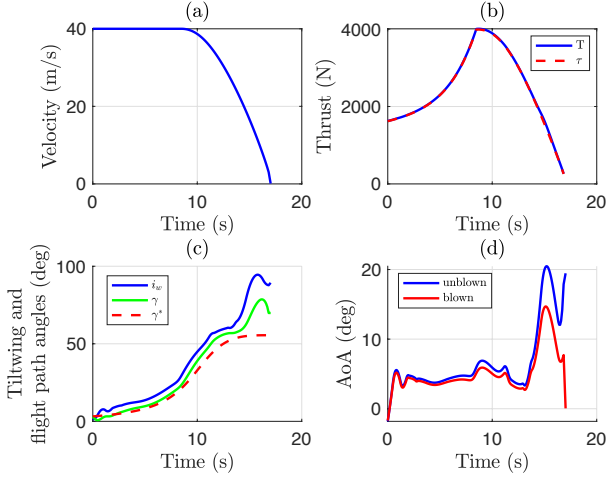} 
         \caption{Backward transition (scenario 2).}
         \label{fig:traj2}
\end{figure}

For completeness, we  consider a second \mcedit{scenario} consisting of a backward transition with an increase in altitude (Figure \ref{fig:traj2}). 
\mcedit{This}
is characterised by \mcedit{an initial}
decrease in velocity magnitude and increase in thrust.
An increase in altitude of about 200\,m is  needed for this manoeuvre due to strict bounds on the effective angle of attack.
A backward transition at constant altitude would require \mcedit{stalling} 
the wing, which is prohibited in the present formulation, illustrating a limitation of our approach. 
To achieve the backward transition, 
a high-drag device or flaps \mcedit{are needed to provide braking forces}.
This was modelled by adding a constant term $d$ to $c(\gamma^{*\prime})$  in problem $\mathcal{P}_1^\dagger$ for the backward transition.

\begin{table}[ht]
\centering
\begin{tabular}{llll}
\hline
\textbf{Parameter} & \textbf{Symbol} & \textbf{Value} & \textbf{Units} \\ \hline
Mass     &    $m$    &   $752.2$    &    \si{kg}   \\ \hline
Gravity acceleration     &    $g$    &   $9.81$    &    \si{m.s^{-2}}   \\ \hline
Wing area &   $S$     &    $8.93$   &    \si{m^2}   \\ \hline
Disk area &   $A$     &    $2.83$   &    \si{m^2} \\ \hline
Wing inertia &   $J_w$     &    $1100$   &    \si{kg.m^2}   \\ \hline
Density of air &   $\rho$     &    $1.225$   &    \si{kg.m^{-3}} \\ \hline
{Lift coefficients} &   $b_0, b_1$   &    $0.43, 0.11$   &    \si{-}, \si{deg^{-1}}  \\ \hline
\multirow{3}*{Drag coefficients} &   $a_0$   &   $0.02$   &    \si{-}  \\
&   $a_1$   &   $0.004$   &    \si{deg^{-1}} \\
&   $a_2$   &   $7.6\mathrm{e}{-5}$   &    \si{deg^{-2}}\\ \hline
Maximum thrust &   $\overline{T}$     &    $8855$   &    \si{N} \\ \hline
Angle of attack range &   $\left[\underline{\alpha}, \overline{\alpha}\right]$     &    $\left[-90, 90\right]$   &    \si{deg}   \\ \hline
Flight path angle range &   $\left[\underline{\gamma}, \overline{\gamma}\right]$     &    $\left[-90, 90\right]$   &    \si{deg}   \\ \hline
Tiltwing angle range &   $\left[\underline{i}_w, \overline{i}_w\right]$     &    $\left[0, 100\right]$   &    \si{deg}   \\ \hline
Acceleration range &   $\left[\underline{a}, \overline{a}\right]$     &    $\left[-0.3g, 0.3g\right]$   &    \si{m . s^{-2}}   \\ \hline
Velocity range &   $\left[\underline{V}, \overline{V}\right]$     &    $\left[0, 40\right]$   &    \si{m/s}   \\ \hline
Momentum range &   $\left[\underline{M}, \overline{M}\right]$     &    $\left[-50; 50\right]$   &    \si{N.m}   \\ \hline
Number of propellers &   $n$     &    $4$   &    \si{-} \\ \hline
Discretisation points &   $N$     &    $1000$   &    \si{-} \\ \hline 
Time step &   $\delta$     &    $0.5$   &    \si{s} \\ \hline 
Degree of polynomial $f$ &   $2n$     &    $26$   &    \si{-} \\ \hline 
\end{tabular}
\vspace{1mm}\caption{Model parameters \mcedit{derived} from A$^3$ Vahana}
\label{tab:param}
\vspace{-3mm}
\end{table}

\begin{table}[ht]
\centering
\begin{tabular}{llll}
\hline
\textbf{Parameter} & \textbf{Symbol} & \textbf{Value} & \textbf{Units} \\ \hline
\multicolumn{4}{c}{{\cellcolor[rgb]{0.753,0.753,0.753}}\textbf{Forward transition}}    \\
Velocity  &   $\left\{V_0; V_f\right\}$     &    $\left\{0.5; 40\right\}$   &    \si{m/s}   \\ \hline
Tiltwing angle  &   $i_0$     &    $75$   &    \si{deg}   \\ \hline
Tiltwing angle rate  &   $\Omega_0$     &    $0$   &    \si{deg/s}   \\ \hline
Flight path angle &   $\gamma_0$     &    $0$   &    \si{deg}   \\ \hline
\multicolumn{4}{c}{{\cellcolor[rgb]{0.753,0.753,0.753}}\textbf{Backward transition}}    \\
Velocity  &   $\left\{V_0; V_f\right\}$     &    $\left\{40; 0.1\right\}$   &    \si{m/s}   \\ \hline
Tiltwing angle  &   $\left\{i_0; i_f\right\}$     &    $\left\{0; 90\right\}$    &    \si{deg}   \\ \hline
Tiltwing angle rate  &   $\Omega_0$     &    $0$   &    \si{deg/s}   \\ \hline
Flight path angle &   $\gamma_0$     &    $1.6$   &    \si{deg}   \\ \hline
\end{tabular}
\vspace{1mm}\caption{Boundary conditions for transitions}
\label{tab:BC}
\vspace{-3mm}
\end{table}

\section{Conclusions}
\label{sec:conclusion}
This paper addresses the trajectory optimisation problem for the transition of a tiltwing VTOL aircraft, leveraging DC decomposition of the dynamics and robust tube programming. The approach is based on successive linearisation of the dynamics around feasible trajectories and treating the linearisation error as a bounded disturbance. The DC form of the dynamics allows to enforce tight bounds on the disturbance via a set of convex inequalities that form the basis of a computationally tractable robust optimisation. The algorithm can compute safe trajectories that are robust to model uncertainty for abrupt transitions at near constant altitude, extending the results in \cite{ACC}. Another contribution of the present work is the  extension of the robust tube optimisation paradigm presented in \cite{DC-TMPC} to dynamic systems that are not convex, by means of a DC decomposition of the nonlinear dynamics. Limitations of the present approach are: i) to obtain a computationally tractable formulation, quadratic approximations of the DC polynomials are required; ii) the computation time, although relatively low compared to solving a NLP, is still too high to leverage the optimisation in a MPC setting. 

Future work will alleviate these problems by i) considering other types of basis functions for the nonlinear dynamics approximation, e.g. radial basis functions that have better scalability than a monomial basis; ii) the use of first order solvers such as ADMM to accelerate computations \cite{me2}. We will then investigate robust MPC for the transition of tiltwing VTOL aircraft. 

\bibliography{biblio} 
\bibliographystyle{ieeetr}

\end{document}